\documentclass[11pt,a4paper]{amsart}
\usepackage[margin=2.7cm]{geometry}
\usepackage[T1]{fontenc}
\usepackage[utf8]{inputenc}
\usepackage{amsmath,amssymb,amsthm}
\usepackage{enumitem}
\usepackage{xcolor}
\usepackage[colorlinks=true,linkcolor=blue!50!black,citecolor=blue!50!black,urlcolor=blue!50!black]{hyperref}

\theoremstyle{plain}
\newtheorem{theorem}{Theorem}
\newtheorem{lemma}{Lemma}
\newtheorem{corollary}{Corollary}
\newtheorem*{maintheorem}{Main Theorem}
\newtheorem*{conjecture}{Conjecture}
\theoremstyle{definition}
\newtheorem{definition}{Definition}
\theoremstyle{remark}
\newtheorem{remark}{Remark}

\newcommand{\Nout}{N^{+}}
\newcommand{\Nin}{N^{-}}
\newcommand{\dout}{d^{+}}
\newcommand{\din}{d^{-}}
\newcommand{\PP}{P_4}
\newcommand{\cp}[5]{#1\leftarrow #2\to #3\leftarrow #4\to #5}

\title{Antidirected paths on five vertices in $5$-chromatic digraphs}
\author{Hussein Karneeb}
\thanks{Lebanese University, KALMA Laboratory, Beirut, Lebanon. Email: hussien.karneeb@st.ul.edu.lb}

\begin{document}

\begin{abstract}
We prove that every $5$-chromatic digraph containing no copy of the regular tournament on five vertices contains both antidirected paths on five vertices. This answers a question posed by El Sahili in $2004$. The proof uses an edge bound for 5-critical graphs to constrain the out-degrees in a counterexample and obtain a contradiction by constructing a proper $4$-colouring.
\end{abstract}
\maketitle

\section{Introduction}\label{sec:intro}

An \emph{oriented path} is an orientation of a path; its \emph{order} is its number of vertices and its \emph{length} is its number of arcs. A \emph{block} is a maximal directed subpath, and $P(k_1,\dots,k_t)$ denotes the oriented path with $t$ blocks of lengths $k_1,\dots,k_t$, starting with a forward block. An antidirected path is an oriented path in which every block has length $1$. On five vertices there are two such paths up to isomorphism $P(1,1,1,1)$ and its converse. A digraph $H$ is \emph{$n$-universal} if every $n$-chromatic digraph contains $H$. The Gallai--Hasse--Roy--Vitaver theorem \cite{Ga68,Ha65,Ro67,Vi62} states that the directed path of order $n$ is $n$-universal. For arbitrary orientations, Burr \cite{Bu80} proved that every $(n-1)^2$-chromatic digraph contains every oriented tree of order $n$ and conjectured that $2n-2$ colours suffice. This bound was improved to $\frac{n^2}{2}-\frac{n}{2}+1$ by Addario-Berry, Havet, Linhares Sales, Reed and Thomass\'e \cite{ABHLRT}, and then to $8\sqrt{2/15}\,n^{3/2}+O(n)$ by Bessy, Gon\c{c}alves and Reinald \cite{BGR}, which is the best published bound. The best bound currently known is a recent arXiv preprint of Fan, Lu and Chen \cite{FLC}, who show that $\lfloor 31\log_2(n!)\rfloor=O(n\log n)$ colours suffice.

Tournaments of order $n$ are $n$-chromatic, and the problem has been completely solved for them. After the work of Gr\"unbaum \cite{Gr71}, Rosenfeld \cite{Ro72} and Thomason \cite{Th86}, Havet and Thomass\'e \cite{HT00} proved Rosenfeld's conjecture: every tournament of order $n\ge 8$ contains every oriented path of order $n$. Three small exceptions exist, namely the directed $3$-cycle, the regular tournament on $5$ vertices and the Paley tournament on $7$ vertices, each missing the antidirected Hamiltonian path. El Sahili \cite{ES15} conjectured that this extends far beyond tournaments.

\begin{conjecture}[El Sahili]
For every $n\ge 8$, every oriented path of order $n$ is $n$-universal.
\end{conjecture}

The conjecture has been studied according to the number of blocks. For paths with two blocks, El Sahili and Kouider \cite{EK07} showed that $(n+1)$ colours suffice. Building on their work, Addario-Berry, Havet and Thomass\'e \cite{AHT07} settled the case of two blocks, proving that every path of order $n\ge 4$ with two blocks is $n$-universal; a new elementary proof was later given by El Sahili, Mortada and Nasser \cite{ENM24}. For three blocks the conjecture remains open: El Joubbeh \cite{EJ22} showed that $(4n-4)$ colours suffice, and a sharper bound of $(n+1)$ colours is known when two consecutive blocks have length $1$ \cite{MEE21}. For four blocks, El Joubbeh and Ghazal \cite{EJG24} obtained the bound $4n-2$, and more generally they proved by induction on $t$ that every oriented path of order $n$ with $t\ge 4$ blocks is contained in every digraph of chromatic number at least $4^{r(t)}n+q(t)$, where $r(t)=\lfloor\log_2(t-1)\rfloor$ and $q(t)=-\frac{20}{7}\,8^{\,r(t)-1}+\frac67$; this improves the bound of \cite{ABHLRT} for paths whose number of blocks is small compared with $n$. Subsequently, Bessy, Gon\c{c}alves and Reinald \cite{BGR} proved that every oriented path of order $n$ with $t\ge 2$ blocks is $\bigl((t-1)(n-3)+3\bigr)$-universal; in particular $2n-3$ colours suffice for three blocks.

The regular tournament $T_5$ on five vertices contains neither antidirected path of order five. El Sahili \cite{ES04} proved that every connected $5$-chromatic oriented graph distinct from $T_5$ with minimum out-degree at least two contains the backward-first antidirected path $P_4$ defined in Section ~\ref{sec:prelim}. He then asked whether a $5$-chromatic oriented graph can contain neither a tournament on five vertices nor $P_4$. We give a negative answer under the weaker assumption that the graph contains no copy of $T_5$. Applying the result to the converse gives the other antidirected path.

\section{Preliminaries}\label{sec:prelim}

All digraphs considered here are finite, without loops, multiple arcs or directed $2$-cycles. The chromatic number of a digraph is that of its underlying graph. Throughout, containment means containment as a subdigraph. For a vertex $x$ we write $\Nout(x)$ and $\Nin(x)$ for its out- and in-neighbourhoods, $\dout(x)=|\Nout(x)|$, $\din(x)=|\Nin(x)|$ and $d(x)=\dout(x)+\din(x)$.

The \emph{antidirected paths of order $5$} are $\PP$, the path $x_1x_2x_3x_4x_5$ oriented as
\[
x_1\leftarrow x_2\to x_3\leftarrow x_4\to x_5 ,
\]
and its converse $P(1,1,1,1)$.

Let $T_5$ be the regular tournament on $5$ vertices. It is $5$-chromatic and contains no $\PP$. We can now state our result precisely.

\begin{maintheorem}
Every $5$-chromatic digraph containing no copy of $T_5$ contains $\PP$.
\end{maintheorem}

Since $T_5$ is isomorphic to its converse, applying the Main Theorem to the converse of a digraph $D$ shows that $D$ also contains the converse of $\PP$, which is the path $P(1,1,1,1)$ of the introduction. Hence $D$ contains both antidirected paths of order $5$.

We use the following results.

\begin{theorem}[El Sahili \cite{ES04}]\label{thm:ES}
Let $D$ be a connected $5$-chromatic digraph, $D\ne T_5$, in which every vertex has out-degree at least $2$. Then $D$ contains $\PP$.
\end{theorem}

A graph $G$ is \emph{$k$-critical} if $\chi(G)=k$ and $\chi(G')<k$ for every proper subgraph $G'$ of $G$.

\begin{theorem}[Kostochka - Yancey \cite{KY14}]\label{thm:KY}
Every $5$-critical graph with $n$ vertices and $m$ edges satisfies $4m\ge 9n-5$.
\end{theorem}

Suppose that the Main Theorem fails for some digraph $D_0$. Its underlying graph contains a $5$-critical subgraph; let $D$ be the corresponding subdigraph of $D_0$. Then $D$ contains neither $\PP$ nor $T_5$. A $5$-critical graph is connected and has minimum degree at least $4$. Writing $n=|V(D)|$ and $E=E(D)$, we therefore have $9n-4|E|\le 5$ by Theorem~\ref{thm:KY}. Throughout Sections~\ref{sec:local}--\ref{sec:proof}, $D$ is this digraph, and we derive a contradiction.

\section{Local structure}\label{sec:local}

Let $H=\{x\in V(D):\dout(x)\ge 3\}$.

\begin{lemma}\label{lem:owner}
Every vertex has at most one in-neighbour in $H$.
\end{lemma}

\begin{proof}
Suppose $x,y\in H$ are distinct in-neighbours of $v$. We may assume $x\not\to y$ after interchanging $x$ and $y$ if necessary, since $D$ has no directed $2$-cycle. Then $|\Nout(x)\setminus\{v,y\}|\ge 2$ and $|\Nout(y)\setminus\{v,x\}|\ge1$, so we can choose $t\in\Nout(y)\setminus\{v,x\}$ and $z\in\Nout(x)\setminus\{v,y,t\}$. Now $\cp zxvyt$ is a copy of $P_4$, a contradiction.
\end{proof}

If $x\in H$ and $x\to z$, we call $x$ the \emph{owner} of $z$; by Lemma~\ref{lem:owner} it is unique.

\begin{lemma}\label{lem:S}
Let $x\in H$, $z\in\Nout(x)$, and $u\in\Nin(z)\setminus\{x\}$. Then either $\Nout(u)\subseteq\{z,x\}$, or $\dout(x)=3$, $x\to u$ and $\Nout(u)=\Nout(x)\setminus\{u\}$.
\end{lemma}

\begin{proof}
Suppose $u\to w$ for some $w\notin\{z,x\}$. If some $t\in\Nout(x)\setminus\{z,u\}$ differs from $w$, then $\cp txzuw$. Hence $\Nout(x)\setminus\{z,u\}=\{w\}$, and since $\dout(x)\ge3$ we get $\Nout(x)=\{z,u,w\}$; in particular $\dout(x)=3$ and $x\to u$. Suppose $u\to t$ for some $t\notin\{z,w\}$. Since $x\to u$, we have $t\neq x$. Then $\cp tuwxz$ is a copy of $P_4$. Hence $\Nout(u)=\{z,w\}=\Nout(x)\setminus\{u\}$.
\end{proof}

\begin{definition}\label{def:special}
A vertex $v\in H$ is \emph{special} if $\dout(v)=3$ and some $u\in\Nout(v)$ satisfies $\Nout(u)=\Nout(v)\setminus\{u\}$; such a $u$ is called a \emph{twin} of $v$, and the set $B(v)=\{v\}\cup\Nout(v)$ is called the \emph{block} of $v$.
\end{definition}

\begin{remark}\label{rem:twin}
Let $v$ be special, $u$ a twin of $v$ and $\Nout(v)=\{u,b,c\}$. Then $u$ is the unique twin of $v$, and $\Nout(u)=\{b,c\}$; in particular $u\notin H$. Indeed, if $u'\ne u$ were another twin, then $u'\in\Nout(v)\setminus\{u\}=\Nout(u)$ and likewise $u\in\Nout(u')$, a $2$-cycle. The rest follows from the definition.
\end{remark}

\begin{definition}\label{def:private}
A vertex $w$ is a \emph{private} of a vertex $z$ if $\Nout(w)=\{z\}$.\\
Let $x\in H$ be the owner of $z\notin H$. We say $z$ is \emph{clean} if every in-neighbour of $z$ other than $x$ is a private of $z$.
\end{definition}

\begin{lemma}\label{lem:nonspecial}
If $x\in H$ is not special, then at most one out-neighbour of $x$ outside $H$ is not clean.
\end{lemma}

\begin{proof}
Let $z\notin H$ be an out-neighbour of $x$ that is not clean, and let $u\in\Nin(z)\setminus\{x\}$ with $\Nout(u)\ne\{z\}$. In Lemma~\ref{lem:S} the second alternative would make $x$ special, so $\Nout(u)\subseteq\{z,x\}$, and hence $\Nout(u)=\{z,x\}$. If $z'\ne z$ is another such vertex, with corresponding $u'$, then $u'\ne u$ (their out-neighbourhoods differ), and $u \neq z'$ otherwise $z'\to x$ , similarly $u' \neq z$, hence $\cp zux{u'}{z'}$ is a copy of $P_4$, a contradiction.
\end{proof}

\begin{lemma}\label{lem:apex}
If a special vertex $v$ has an in-neighbour in $H$, then its twin $u$ is clean.
\end{lemma}

\begin{proof}
Write $\Nout(v)=\{u,b,c\}$. Let $x\in H$ with $x\to v$, and suppose $w\ne v$ is an in-neighbour of $u$ that is not a private of $u$. Note that $w\notin\{b,c\}$, since $u\to b$, $u\to c$ but $w\to u$.

If $w\to t$ with $t\notin\{u,v\}$, choose $t'\in\{b,c\}\setminus\{t\}$; then $\cp twuv{t'}$ is a copy of $P_4$. Hence $\Nout(w)\subseteq\{u,v\}$, and as $w$ is not a private of $u$, $\Nout(w)=\{u,v\}$; in particular $w\notin H$, so $w\ne x$.
Choose $w'\in\Nout(x)\setminus\{v,w\}$, which is possible since $\dout(x)\ge3$. Then $w'\ne u$: otherwise $u$ would have the two in-neighbours $x$ and $v$ in $H$, contradicting Lemma~\ref{lem:owner}. Hence $\cp uwvx{w'}$ is a copy of $\PP$, a contradiction.
\end{proof}

\begin{remark}\label{rem:private}
Every clean vertex has a private, and distinct vertices have distinct privates. Indeed, a clean vertex $z$ is not in $H$, so $\dout(z)\le 2$; since $d(z)\ge4$, we get $\din(z)\ge 2$. Thus $z$ has an in-neighbour other than its owner, which is a private of $z$. The second statement is clear.
\end{remark}

\section{Counting}\label{sec:count}

Let $h$ be the number of special vertices, $X$ the set of non-special vertices of $H$, $L=\{x\in V(D):\dout(x)=1\}$, $S=\{x\in V(D):\dout(x)=0\}$, $R$ the set of vertices lying in no block, and $s$ the number of special vertices having a special in-neighbour.

\begin{lemma}\label{lem:exact}
\[
9n-4|E|=|R|-s+4|L|+8|S|-4\sum_{x\in X}\bigl(\dout(x)-2\bigr).
\]
\end{lemma}

\begin{proof}
Since $|E|=\sum_x\dout(x)$, we have $9n-4|E|=n+4\sum_x(2-\dout(x))$. By Lemma~\ref{lem:owner}, the out-neighbourhoods of distinct special vertices are disjoint, so their union $W$ has $3h$ elements; and a special vertex lies in $W$ if and only if it has a special in-neighbour. Hence the blocks cover exactly $h+3h-s$ vertices, and $n=4h-s+|R|$.
Splitting $\sum_x(2-\dout(x))$ according to whether $x$ is special, in $X$, of out-degree $2$, in $L$ or in $S$ gives
\[
\sum_x\bigl(2-\dout(x)\bigr)=-h-\sum_{x\in X}\bigl(\dout(x)-2\bigr)+|L|+2|S|,
\]
and the claim follows.
\end{proof}

Let $Z$ be the set of vertices outside $H$ that are owned by a vertex of $X$ and are not clean, and let $X_1$ be the set of vertices of $X$ having an in-neighbour in $X$. By Lemmas~\ref{lem:owner} and~\ref{lem:nonspecial}, the map sending $z\in Z$ to its owner is an injection $Z\to X$; in particular
\begin{equation}\label{eq:ZX}
|X_1|\le|X|\quad\text{and}\quad|Z|\le|X|.
\end{equation}

For $x\in H$ we define a set $Y(x)$ of vertices, in two cases.
\begin{enumerate}[label=\textup{(\roman*)},nosep]
\item If $x$ is special, $Y(x)$ is the set of twins of the special out-neighbours of $x$.
\item If $x\in X$, $Y(x)$ is the set of twins of the special out-neighbours of $x$, together with the set of clean out-neighbours of $x$ outside $H$.
\end{enumerate}
Let $P(x)$ be the set of vertices $w$ with $\Nout(w)=\{z\}$ for some $z\in Y(x)$, that is, the set of privates of the vertices in $Y(x)$.

\begin{lemma}\label{lem:Pdisjoint}
The sets $P(x)$, $x\in H$, are pairwise disjoint subsets of $L$, and $|P(x)|\ge|Y(x)|$ for every $x\in H$.
\end{lemma}

\begin{proof}
Each $w\in P(x)$ has out-degree $1$, so $P(x)\subseteq L$, and $w$ determines the vertex $z\in Y(x)$ with $\Nout(w)=\{z\}$. We claim that $z$ determines $x$. If $z$ is the twin of a special out-neighbour $v$ of $x$, then $v$ is the owner of $z$ and $x$ is the owner of $v$. Otherwise $z$ is a clean out-neighbour of $x\in X$ outside $H$, so $x$ is itself the owner of $z$, and it is not special. By Lemma~\ref{lem:owner} these owners are unique, and the two situations cannot both occur for the same $z$, since the owner of $z$ is special in the first and not special in the second. This proves the disjointness.

For the inequality, every $z\in Y(x)$ is clean: by definition if $z$ is a clean out-neighbour of $x$ outside $H$, and by Lemma~\ref{lem:apex} if $z$ is the twin of a special out-neighbour $v$ of $x$, since $v$ has the in-neighbour $x\in H$. Hence, by Remark~\ref{rem:private}, every $z\in Y(x)$ has a private. Since distinct vertices have distinct privates, $|P(x)|\ge|Y(x)|$.
\end{proof}

\begin{lemma}\label{lem:L}
$\displaystyle |L|\ \ge\ s+\sum_{x\in X}\dout(x)-|X_1|-|Z|$.
\end{lemma}

\begin{proof}
By Lemma~\ref{lem:Pdisjoint},
\[
|L|\ \ge\ \sum_{x\in H}|P(x)|\ \ge\ \sum_{x\in H}|Y(x)|\ =\ \sum_{x\text{ special}}|Y(x)|\ +\ \sum_{x\in X}|Y(x)|.
\]

For \emph{special vertices:} Distinct special vertices have disjoint out-neighbourhoods (Lemma~\ref{lem:owner}), hence distinct twins, so $|Y(x)|$ is the number of special out-neighbours of $x$. Summing over special $x$ counts the arcs between special vertices. By Lemma~\ref{lem:owner} each special vertex is the head of at most one such arc, so the number of arcs equals the number of special vertices having a special in-neighbour. Thus $\sum_{x\text{ special}}|Y(x)|=s$.

For \emph{vertices of $X$:} Let $x\in X$. Every out-neighbour of $x$ is exactly one of the following: special; in $X$; outside $H$ and clean; or outside $H$ and not clean, that is in $Z$. The twins of the special out-neighbours of $x$ are not out-neighbours of $x$ (otherwise a twin would have two in-neighbours in $H$, namely $x$ and its owner $v\ne x$). Hence the two parts of $Y(x)$ are disjoint, and
\[
|Y(x)|\ =\ \dout(x)-|\Nout(x)\cap X|-|\Nout(x)\cap Z|.
\]
By Lemma~\ref{lem:owner}, every vertex of $X_1$ has exactly one in-neighbour in $X$, and every vertex of $Z$ has exactly one owner, which lies in $X$. Therefore $\sum_{x\in X}|\Nout(x)\cap X|=|X_1|$ and $\sum_{x\in X}|\Nout(x)\cap Z|=|Z|$, and summing over $x\in X$ gives $\sum_{x\in X}|Y(x)|=\sum_{x\in X}\dout(x)-|X_1|-|Z|$.

Combining the two sums, we obtain
\[
|L|\ \ge\ \sum_{x\text{ special}}|Y(x)|+\sum_{x\in X}|Y(x)|\ \ge\ s+\sum_{x\in X}\dout(x)-|X_1|-|Z|,
\]
which proves the lemma.
\end{proof}

\begin{corollary}\label{lem:key}
\[
5\ \ge\ 9n-4|E|\ \ge\ |R|+3s+8|S|+4\bigl(|X|-|X_1|\bigr)+4\bigl(|X|-|Z|\bigr),
\]
and all terms on the right-hand side are non-negative.
\end{corollary}

\begin{corollary}\label{cor:Ss}
$S=\emptyset$ and $s\le1$.
\end{corollary}

\begin{proof}
By Corollary~\ref{lem:key}, all terms on its right-hand side are non-negative, so $5\ge 8|S|$ and $5\ge 3s$. Hence $|S|=0$ and $s\le 1$.
\end{proof}

\section{Proof of the Main Theorem}\label{sec:proof}

We will prove that $X=\emptyset$.

\begin{lemma}\label{lem:block}
Let $v$ be special with twin $u$ and $\Nout(v)=\{u,b,c\}$.
\begin{enumerate}[label=\textup{(\alph*)}]
\item Every in-neighbour of $b$ other than $v$ and $u$ is a private of $b$; similarly for $c$.
\item If $b\in X$, then every out-neighbour of $b$ outside $H$ is clean; similarly for $c$.
\end{enumerate}
\end{lemma}

\begin{proof}
(a) Let $w\in\Nin(b)\setminus\{u,v\}$ and suppose $w\to t$ with $t\ne b$; we find a copy of $\PP$ in each of the following cases:
\begin{itemize}[nosep]
\item If $w=c$, then $c\to t$ and the arcs $u\to c$ and $v\to c$ imply $t\notin \{u,v\}$. Since $t\neq b$ and $t\neq c$, the path $\cp tcbvu$ is a copy of $P_4$.
\item If $w\ne c$, then $w\notin\{b,c,u,v\}$, and we distinguish three cases according to $t$. If $t\notin\{u,v\}$, then $\cp twbvu$. If $t=u$, then $\cp bwuvc$. If $t=v$, then $\cp vwbuc$.
\end{itemize}
In every case we reach a contradiction, so $w$ is a private of $b$.

(b) Note that $b\not\to c$ by Lemma~\ref{lem:owner} (since $c$ already has the in-neighbour $v\in H$, and $b\in H$, $b\ne v$). Let $t\notin H$ be an out-neighbour of $b$, so that $b$ is the owner of $t$, and suppose some $w\in\Nin(t)\setminus\{b\}$ is not a private of $t$. As $b$ is not special, the second alternative of Lemma~\ref{lem:S} (applied to $x=b$, $z=t$ and $u=w$) cannot occur. Hence $\Nout(w)\subseteq\{t,b\}$, and since $w$ is not a private of $t$, $\Nout(w)=\{t,b\}$. In particular $\dout(w)=2$ and $w\to b$. By part~(a), $w\in\{u,v\}$ or $w$ is a private of $b$; the latter is impossible as $\dout(w)=2$, and $w=v$ is impossible as $\dout(v)=3$. So $w=u$. But then $t\in\Nout(u)=\{b,c\}$, which is impossible since $t\ne b$ and $b\not\to c$.
\end{proof}

\begin{lemma}\label{lem:X}
$X=\emptyset$.
\end{lemma}

\begin{proof}
Suppose that $X\ne\emptyset$. Suppose first that some $x\in X$ lies in a block $B(v)$. Then $x$ is neither $v$ (which is special) nor the twin of $v$ (which is not in $H$), so $x$ is a non-twin out-neighbour of $v$. Its only in-neighbour in $H$ is $v$, which is not in $X$, so $x\notin X_1$; and by Lemma~\ref{lem:block}(b), $x$ owns no vertex of $Z$, so $x$ is not in the image of the injection $Z\to X$. Thus $|X|-|X_1|\ge1$ and $|X|-|Z|\ge1$, and Corollary~\ref{lem:key} (all of whose terms are non-negative) gives $5\ge 4+4=8$, which is impossible. Hence $X\subseteq R$.

Also $Z\subseteq R$: a vertex $z\in Z$ in a block $B(v)$ would be an out-neighbour of $v$ (as $z\notin H$), so its owner would be the special vertex $v$, not a vertex of $X$.

Since $X\subseteq H$ and $Z\cap H=\emptyset$, we get $|R|\ge|X|+|Z|$, and Corollary~\ref{lem:key} together with $|Z|\le|X|$ yields
\[
5\ \ge\ |X|+|Z|+8|X|-4|X_1|-4|Z|\ =\ 9|X|-4|X_1|-3|Z|\ \ge\ 6|X|-4|X_1|.
\]
We show that $6|X|-4|X_1|\ge 6$, a contradiction.
\begin{itemize}[nosep]
\item If $X_1\ne X$, then $|X_1|\le|X|-1$, so $6|X|-4|X_1|\ge 2|X|+4\ge 6$, because $|X|\ge1$.
\item If $X_1=X\ne\emptyset$, then every vertex of $X$ has an in-neighbour in $X$. As $D$ has no loops and no $2$-cycles, $|X|\ge3$. Then $6|X|-4|X_1|=2|X|\ge 6$.
\end{itemize}
Both cases give $5\ge6$, which is impossible, so $X=\emptyset$.
\end{proof}

\begin{proof}[\textbf{Proof of the Main Theorem}]
By Lemma~\ref{lem:exact}, Corollary~\ref{cor:Ss} and Lemma~\ref{lem:X},
\begin{equation}\label{eq:final}
5\ \ge\ 9n-4|E|\ =\ |R|-s+4|L|,\qquad s\le1 .
\end{equation}
Hence $4|L|\le 5+s-|R|\le 6$, so $|L|\le 1$. If $L=\emptyset$, then, as $S=\emptyset$ too, no vertex of $D$ has out-degree $0$ or $1$, i.e.\ every vertex of $D$ has out-degree at least $2$, and Theorem~\ref{thm:ES} shows that $D$ contains $\PP$ (recall that $D$ is connected and $D\ne T_5$), a contradiction. So $L=\{p\}$ for some vertex $p$, and~\eqref{eq:final} gives $5\ge |R|-s+4$, that is $|R|\le s+1\le2$.

We claim $p\in R$. Suppose not. Then $p$ lies in a block $B(v)$. Since $\dout(p)=1$, while special vertices have out-degree $3$ and twins have out-degree $2$, $p$ is a non-twin out-neighbour of the special vertex $v$; write $\Nout(v)=\{u,b,c\}$ with $u$ the twin of $v$, and say $p=b$ (the case $p=c$ is symmetric). The minimum degree of $D$ is at least $4$, so $\din(p)\ge 4-\dout(p)=3$. Two of the in-neighbours of $p$ are $v$ and $u$, so $p$ has a third in-neighbour $w\notin\{u,v\}$. By Lemma~\ref{lem:block}(a), $w$ is a private of $p$. Then $w\in L$ and $w\ne p$, contradicting $|L|=1$.

We now construct a proper $4$-colouring of $D$. Every vertex of a block is either special, a twin, or neither. We colour the special vertex having a special in-neighbour, if it exists, with colour $4$ (by Corollary~\ref{cor:Ss} there is at most one), all other special vertices with colour $1$, all twins with colour $2$, and all remaining block vertices with colour $3$.

Let $x\to y$ be an arc between block vertices of the same colour. The colour is not $1$, since then $y$ would have colour $4$, nor $4$, since $s\le1$. If it is $2$, then $x$ and $y$ are twins of special vertices $v$ and $v'$; since $\Nout(x)=\Nout(v)\setminus\{x\}$, we get $v\to y$, so $v=v'$ by Lemma~\ref{lem:owner}, and $x=y$ by uniqueness of the twin, which is impossible. If it is $3$, then $y$ is a non-twin out-neighbour of a special vertex $v$, and $x$ is neither $v$ nor its twin; by Lemma~\ref{lem:block}(a), $x$ is a private of $y$, so $x\in L=\{p\}$, contradicting $p\in R$. Thus the colouring of the blocks is proper.

The out-neighbours of a special vertex lie in its block, and so do those of a twin. Hence every arc from a block vertex to a vertex of $R$ starts at a vertex of colour $3$. If $R\setminus\{p\}$ is nonempty, let $r$ be its unique vertex and colour $r$ first. Since $X=\emptyset$ and every special vertex belongs to a block, $r\notin H$, so $d^+(r)\le 2$. Its already coloured in-neighbours all have colour $3$, and its out-neighbours forbid at most two further colours. Thus a colour is available for $r$. Finally, the neighbours of $p$ are its in-neighbours in blocks (colour $3$), its unique out-neighbour and possibly $r$, so again some colour is available for $p$. This yields a proper $4$-colouring of $D$, the desired contradiction.
\end{proof}

\end{document}